\documentclass[12pt,a4paper]{article}
\usepackage{amssymb}
\usepackage{amsfonts}
\usepackage{srcltx}
\usepackage{amsmath}
\usepackage{amsthm}
\usepackage{enumerate}

\newtheorem{theorem}{Theorem}[section]
\newtheorem{proposition}[theorem]{Proposition}
\newtheorem{lemma}[theorem]{Lemma}

\newtheorem{remark}{Remark}[section]

\newcommand\cA{{\cal A}}
\newcommand\cC{{\cal C}}

\newcommand\cF{{\cal F}}
\newcommand\cL{{\cal L}}
\newcommand\cB{{\cal B}}

\newcommand\cN{{\cal N}}

\newcommand\cX{{\cal X}}

\newcommand\cR{{\cal R}}

\newcommand\cU{{\cal U}}

\newcommand\e{\epsilon}
\newcommand\ve{\varepsilon}
\newcommand\ov{\overline}

\def\bbr{{\mathbb R}}
\def\bbn{{\mathbb N}}

\def\text#1{\hbox{#1}}

\def\endproof{\mbox{\ $\qed$}}

\def\E{{\bf E}}
\def\P{{\bf P}}
\def\C{{\bf C}}

\def\H{{\bf H}}

\def\M{{\bf M}}

\def\Chi{{\bf 1}}

\def\d{\mathrm{d}}
\def\build #1_#2{\mathrel{\mathop{\kern 0pt #1}\limits_{#2}}} 
\newcommand{\zs}[1]{{\mathchoice{#1}{#1}{\lower.25ex\hbox{$\scriptstyle#1$}}
{\lower0.25ex\hbox{$\scriptscriptstyle#1$}}}}

\numberwithin{equation}{section}
\begin{document}
\title{
Adaptive sequential estimation for
ergodic diffusion processes in quadratic metric.\\
{\Large Part 2: Asymptotic efficiency.}
}
\author{{\Large By  L.Galtchouk and S.Pergamenshchikov}
\thanks{The second author is partially supported by the RFFI-Grant 04-01-00855.}\\
University  of Strasbourg and University of Rouen
}
\date{}
\maketitle

\begin{abstract}

Asymptotic efficiency is proved for the constructed in part 1
procedure, i.e. Pinsker's
constant is found in the asymptotic lower bound for the minimax quadratic risk. It is shown that
the asymptotic minimax quadratic risk of the constructed procedure coincides with this
constant. 
\end{abstract}
\bibliographystyle{plain}

\newpage

\section{Introduction}\label{Se.1}

%\Pe\marginpar{\Pe} 
The paper is a continuation of the investigation carried in \cite{GaPe6}
and it deals with asymptotic nonparametric estimation of the drift coefficient
 $S$ in observed diffusion process $(y_t)_{t\ge 0}$ governed by the stochastic differential
 equation 
%We consider the problem of estimating the drift coefficient of diffusion process
 \begin{equation}\label{1.1}
\d y_t=S(y_t)\,\d t +\d w_t\,,\quad 0\le t\le T\,,\quad y_0=y\,,
 \end{equation}
where $(w_t)_{t\ge 0}$ is a scalar 
standard Wiener process, $y_0=y$ is a initial condition.
% and  the function $S(\cdot)$ is unknown function. 
%The problem is to estimate the function $S$ from observations of $(y_t)_\zs{0\le t\le T}$.

In the paper \cite{GaPe6} we have constructed a non-asymptotic adaptive procedure
for which a  sharp non-asymptotic oracle inequality is obtained. This oracle inequality gives
a upper bound for a quadratic risk. In this paper we analyze asymptotic properties  (as $T\to\infty$)
of the above adaptive procedure and state that it is asymptotically efficient. This means that the 
procedure provides the optimal convergence rate and the best constant (the Pinsker constant).
 
The problem of asymptotic (as $T\to\infty$) minimax nonparametric estimation
of the drift coefficient $S$ in the model \eqref{1.1} has been studied 
in a number of papers, see for example, \cite{Da}--\cite{GaPe5}. So the papers
 \cite{GaPe1}, \cite{GaPe3} and \cite{GaPe5} deal with the estimation problem at a fixed point.
 In \cite{GaPe3} and \cite{GaPe5} in the case of known smoothness of the function $S$,  efficient procedures were
 constructed which possess the optimal convergence rate and which provide the sharp minimax constant
 in asymptotic risks. Further in \cite{GaPe3}, a adaptive estimation procedure was given when the
   smoothness of the function $S$ is unknown, the procedure provides the optimal convergence rate.
   Moreover, for estimation in $\cL_2-$norme, in \cite{GaPe2} a adaptive sequential estimation procedure was
   constructed. The procedure possesses the optimal convergence rate and it is based on the model 
   selection (see, \cite{BaBiMa} and \cite{FoPe}).

%{\Large\Pe\marginpar{\Pe} 
%Tak v rabotah \cite{GaPe1}, \cite{GaPe3} and \cite{GaPe5} byla rassmotrtena zadacha ocenivaniya v tochke.
% V  \cite{GaPe3} and \cite{GaPe5}
%dlya sluchaya izvestnoi gladkosti byli postroeny effectivnye procedures 
%obladayushie optimal'noi skorost'yu shodimosti i dostigayeshei tochnoi minimal'noi konstanty
%v asimptoticheskih riskah. Dalee v 
%\cite{GaPe3} byla postroena adaptivnaya procedura ocenivanya dlya sluchaya neizvestnoi
%regulyarnosti funkcii $S$ dlya kotoroi byla pokazana minimaxnost' po skorosti shodimosti.
%Bolee togo, dlya ocenivaniya v $\cL_2-$ norme  v \cite{GaPe2}
% byla postroena adaptivnaya posledovatel'naya procedura ocenivaniya
%na osnove vybora modeli (see,\cite{BaBiMa} and \cite{FoPe}) dlya
%kotoroi byla dokazana optimal'naya skorost' shodimosti.
%}

%For the drift estimation in a known H\"older class the optimal convergence 
%rate of the minimax risk 
%was obtained in \cite{GaPe3} and \cite{GaPe5}
% for the pointwise
%estimation 
%and in \cite{GaPe2} for the
%$L_2([a,b],dx)-$loss function.
%When the H\"older smoothness of the drift $S$ is unknown, an adaptive
% estimator possessing the optimal convergence rate is proposed in
% \cite{GaPe1} for the pointwise
%estimation. 

The sharp asymptotic bounds and efficient estimators for the drift $S$ in 
model (\ref{1.1}) with the known Sobolev smoothness
 was given in \cite{DaKu} and with unknown one in \cite{Da}
 for local weighted  $\cL_2-$losses, where the weight function is equal to the
squared unknown ergodic density.
Note that the weighted  $\cL_2-$risk considered in the papers \cite{Da}-\cite{DaKu}
is restrictive for the following reasons. The ergodic density being exponentially
decreasing, the feasible estimation is possible on an finite interval which
 depends on unknown function $S$. Moreover, the weighted  $\cL_2-$risk in 
these papers is local and the centres of vicinities in the local risk should be
smoother than the function to be estimated.
Since in the local risk the vicinity radius tends to zero, it means really that the proposed
 procedure estimates the 
centre of the vicinity which can be estimated with a better convergence rate.
% We thing that 
So the approach proposed in \cite{Da}-\cite{DaKu} 
%has not any statistic importance because it
permets to calculate the sharp asymptotic constant by lossing the optimal convergence rate.

%{\Large\Pe\marginpar{\Pe} 
%Tak kak v lokal'nom riske raduis okrestnosti stremit'sya k nulu, to 
%eto oznachjaet, chto fakticheski ocenivaetsya centre okrestnosti,
%kotoryi mojet byt' ocenen s luchshei skorostiy.
%Na nash vzglyad podhod predlojennyi v \cite{Da}-\cite{DaKu} 
%kotoryi za schet poteri skorosti shodimosti pozvolyaet podchitat' asymptoticheskuyu
%konstantu
%lishen statisticheskogo smysla.
%}

 In this paper we consider the global $\cL_2-$risk and we show how to obtain the optimal
 convergence rate and to reach the Pinsker constant. We prove that the  constructed
 in \cite{GaPe6} procedure provides the both above properties.

%{\Large\Pe\marginpar{\Pe} 
%V dannoi je rabote my pokazyvaem kak sohranit' optimal'nuyu skorost' shodimosti
%i dostich konstantu Pinskera. A imenno dokazyvaetsya, chto postroennaya 
%v
%rabote \cite{GaPe6} procedura shoditsya s minimaxnoi skorost'yu i pri etom dostigaetsya
%konstanta Pinskera.
%}

%Such a minimax criterion has no statistical sense when the radius of vicinities tends to zero.

%Zametim je , chto na samom dele risk predlojennyi v rabotah \cite{DaKu}
% \cite{Da} yavlaetsya ogranichitel'nym, v tom smysle, chto ergodicheskaya
%plotnost' exponencial'no ubyvaet, a sledovatel'no, real'noe ocenivanie proishodit
%na konechnom intervale, kotoryi k stati zavisit ot neizvestnoi funkcii. 
%Bolee, togo ocenivanie proishodit v lokal'nom smysle, chto toje ogranichivaet
%oblast' primeneniya, tak kak na samom dele
%trebuetsya okrestnost' s central'noi funksiei bolee gladkoi, chem
%ocenivaemaya. Takoi kriterii v asimptoticheskom smysle, kogda radius centra
%stremista k nulu v obshe ne immet nikakogo statisticheskogo smysla.
%}

The paper is organized as follows. 
In the next Section we formulate the problem
and give the definitions of the functional classes
 and the global quadratic risk.
In Section 3 the sequential adaptive procedure is constructed. The sharp upper
 bound for
the global minimax quadratic risk over all estimates is given in Section 4 (Th. 4.1).
 In 
Section 5 we prove that the lower bound of the global risk for 
the sequential kernel 
estimate coincides with the sharp lower bound, i.e. 
this estimate is asymptotically efficient. The Appendix contains the proofs of 
auxiliairy results.

\medskip

\section{Main results.}\label{Se.2}

Let $(\Omega,\cF,(\cF_t)_{t\ge 0}, \P )$ be a filtered 
probability space
satisfying the usual conditions and $(w_t, \cF_t)_{t\ge 0}$ be a standard Wiener
process.

Suppose that the observed process $(y_t)_{t\ge 0}$ is governed by the
stochastic differential equation (\ref{1.1}),
where the unknown function $S(\cdot)$ satisfies the Lipschitz condition,
$S\in Lip_\zs{L}(\bbr)$, with
$$
Lip_\zs{L}(\bbr)\,=\,
\left\{f\in\cC(\bbr)\,:\,\sup_\zs{x,y\in\bbr}\frac{|f(x)-f(y)|}{|x-y|}\,\le\,L\right\}\,.
$$
In this case the
equation \eqref{1.1} admits a strong solution.
 We denote by
$(\cF^y_t)_{t\ge 0}$ the natural filtration of the process
 $(y_t)_{t\ge 0}$ and by $\E_\zs{S}$ the expectation with respect to the distribution 
 law $\P_\zs{S}$ of the process $(y_t)_{t\ge 0}$ given the drift $S$.
The problem is to estimate the function $S$ in $\cL_2[a,b]$-risk, 
for some $a<b, b-a\ge 1$, 
i.e. for any estimate $\hat{S}_\zs{T}$ of $S$ based on 
$(y_t)_\zs{0\le t\le T}$, we consider the following quadratic risk :
\begin{equation}\label{2.1}
\cR(\hat{S}_\zs{T},S)\,=\,\E_\zs{S}\|\hat{S}_\zs{T}-S\|^2\,,\quad
\|S\|^2\,=\,\int^b_\zs{a}\,S^2(x)\,\d x\,.
\end{equation}

To obtain a good  estimate for the function 
$S$ it is
 necessary to impose some conditions on the function 
$S$ which are 
similar to the periodicity of the  deterministic signal
  in the  white noise model.
    One of conditions which is sufficient for this
purpose is the assumption that the process
$(y_t)$ in (\ref{1.1})  returns to any vicinity of each point $x\in [a,b]$
 infinite times. 
The ergodicity of $(y_t)$ provides this property.

Let $L\ge 1$.
We define the following functional class :
 \begin{align}\nonumber
\Sigma_\zs{L}=\{&S\in Lip_\zs{L}(\bbr)\,:\,\sup_\zs{|z|\le L}\,|S(z)|\le L \,; \ 
 \forall |x|\ge L\,, \\ \label{2.2}
& \exists\ \dot{S}(x)\in \cC\ \mbox{\rm such that}\ -L\le \dot{S}(x)\le-1/L\}\,.
 \end{align}
It is easy to see that 
\begin{equation}\label{2.2-1}
%s_\zs{*}\,=\,\sup_\zs{x\in[a,b]}\,
\nu^{*}\,=\,\sup_\zs{x\in[a,b]}\,
\sup_\zs{S\in\Sigma_\zs{L}}\,S^2(x)\,<\,\infty\,.
 \end{equation}
Moreover, if $S\in \Sigma_\zs{L}$, 
then there exists the ergodic density
 \begin{equation}\label{2.3}
q(x)\,=q_\zs{S}(x)\,=\,\frac{\exp\{2\int^x_0 S(z)\d z\}}{
\int^{+\infty}_{-\infty}\exp\{2\int^y_0 S(z)\d z\}\d y}
 \end{equation}
(see,e.g., Gihman and Skorohod (1972), Ch.4, 18, Th2). 
It easy to see that
this density satisfies the following inequalities
\begin{equation}\label{2.4}
0<\,q_*:=\,\inf_\zs{|x|\le b_\zs{*}}\,\inf_\zs{S\in\Sigma_\zs{L}}\,q_\zs{S}(x)\,
\le\,\sup_\zs{x\in\bbr}\,\sup_\zs{S\in\Sigma_\zs{L}}\,q_\zs{S}(x):=
q^*<\infty\,,
\end{equation}
where $b_\zs{*}=1+|a|+|b|$.
Let $S_\zs{0}$ be a known $k$ times differentiable function from $\Sigma_\zs{L}$.
We define the following functional class
\begin{equation}\label{2.5}
W^{k}_\zs{r}=\{S\,:\,S-S_\zs{0}\in \Sigma_\zs{L}\cap\,\cC^{k}_\zs{per}([a,b])
\,,\,\sum_\zs{j=0}^k\,\|S^{(j)}-S^{(j)}_\zs{0}\|^2\le r\}\,,
 \end{equation}
where $r>0\,,\ k\ge 1$ are 
some  parameters, $\cC^{k}_\zs{per}([a,b])$ is a set of 
 $k$ times differentiable functions 
$f\,:\,[a,b]\to\bbr$ such that
$f^{(i)}(a)=f^{(i)}(b)$ for all $0\le i \le k$.

Note that, we can represent the functional class $W^{k}_\zs{r}$ as the ellipse in $\cL_\zs{2}[a,b]$, i.e.
 \begin{equation}\label{2.6}
W^{k}_\zs{r}=\{S\,:\,S-S_\zs{0}\in \Sigma_\zs{L}\cap\,\cC^{k}_\zs{per}([a,b])
\,,\,\sum_\zs{j=1}^{\infty}\,\varpi_\zs{j}\,\theta^2_\zs{j}\,\le r\}\,,
 \end{equation}
where 
$$
\varpi_\zs{j}=\sum^k_{i=0}\left(\frac{2\pi [j/2]}{b-a}\right)^{2i}
\quad\mbox{and}\quad
 \theta_\zs{j}=\int^b_\zs{a}(S(x)-S_\zs{0}(x))\phi_\zs{j}(x)\d x\,.
$$
Here $(\phi_\zs{j})_\zs{j\ge 1}$ is the standard trigonometric basis in $\cL_\zs{2}[a,b]$ 
(see the definition (4.6) in \cite{GaPe6}) and $[a]$ is the integer part of a number $a$.

\begin{remark}
Note that the functional class $W^{k}_\zs{r}$ is an ellipse with the centre at $S_\zs{0}$.
Usually in  such kind problems one takes an ellipse with the centre $S_\zs{0}\equiv 0$. In
the model \eqref{1.1} we cannot take $S_\zs{0}\equiv 0$ as the centre since  this function 
does not belong to the space $\Sigma_\zs{L}$, i.e. the process \eqref{1.1} is not ergodic
for this function.
%{\Large
%\Pe\marginpar{\Pe} 
%Zametim, chto mnojestvo $W^{k}_\zs{r}$ est ellips v centre $S_\zs{0}$. Obuchno b takih zadachah
%beretsya ellips s centrom $S_\zs{0}=0$. No v modele  \eqref{1.1} my ne mojem vzyat'
%sentre $S_\zs{0}=0$ tak kak nul' ne vhodit v $\Sigma_\zs{L}$, t.e. 
%dlya nulevoi funkcii process \eqref{1.1} ne yavlyaetsya ergodicheskim. 
%}
\end{remark}

In \cite{GaPe6} we constructed an adaptive sequential estimator $\hat{S}_\zs{*}$ for which 
 the oracle inequality (Theorem 4.2) holds.  In this paper we prove that this inequality is
sharp in the asymptotic sense, i.e. we show that the minimax
quadratic risk for $\hat{S}_\zs{*}$
asymptotically equals to the Pinsker constant.

 To formulate our asymptotic results
we define the following normalizing coefficient 
\begin{equation}\label{2.7}
\gamma(S)\,=\,((1+2k)r)^{1/(2k+1)}\,
\left(\frac{J(S)k}{\pi (k+1)}\right)^{2k/(2k+1)}
\end{equation}
with
\begin{equation}\label{2.8}
J(S)\,=\,\int_a^b\frac{1}{q_\zs{S}(u)}\,\d u\,.
\end{equation}
It is well known that for any $S\in W^{k}_\zs{r}$ 
 the optimal rate is 
$T^{2k/(2k+1)}$ (see, for example, \cite{GaPe2}).
In this paper 
we show that the adaptive estimator $\hat{S}_\zs{*}$, defined by (4.17)
in \cite{GaPe6}, is asymptotically efficient.
\begin{theorem}\label{Th.2.1} 
The quadratic risk \eqref{2.1} of the
 sequential estimator $\hat{S}_\zs{*}$ has the following asymptotic upper bound
 \begin{equation}\label{2.9}
\limsup_{T\to\infty}\,T^{2k/(2k+1)}
\sup_\zs{S\in W^k_r}\,
\frac{\cR(\hat{S}_\zs{*},S)}{\gamma(S)}
\le 1 \,.
 \end{equation}
\end{theorem}

Moreover, the following result claims that this upper bound is sharp.
\begin{theorem}\label{Th.2.2}
For any estimator $\hat{S}_\zs{T}$ of $S$ measurable with respect to $\cF^y_T$,
 \begin{equation}\label{2.10}
\liminf_{T\to\infty}\inf_{\hat{S}_\zs{T}}\,T^{2k/(2k+1)}\sup_\zs{S\in W^k_r}
\frac{\cR(\hat{S}_\zs{T},S)}{\gamma(S)}
\ge 1 \,.
 \end{equation}
\end{theorem}
Our approach is based on the truncated sequential procedure proposed in
\cite{GaPe1}, \cite{GaPe2} and  \cite{GaPe5} for the diffusion model 
(\ref{1.1}). Through this procedure we pass to discrete regression model
in which we make use of the adaptive procedure $\hat{S}_\zs{*}$   proposed in 
\cite{GaPe4} for the family $(\hat{S}_\zs{\alpha}\,,\alpha\in\cA)$, 
where $\hat{S}_\zs{\alpha}$ is a weighted least squares estimator
with the Pinsker weights. In the next section we describe the discrete regression
model.

\medskip

\section{Adaptive procedure}\label{Se.AD}

We remind of that in \cite{GaPe6} we pass by the sequential method  to discrete scheme
at the points
\begin{equation}\label{AD.1}
x_\zs{l}=a+\frac{l}{n}(b-a)\,, \quad 1\le l\le n\,,
\end{equation}
with $n=2[(T-1)/2]+1$.
At each $x_\zs{l}$ we use the sequential kernel estimator
\begin{equation}\label{AD.1-2}
\left\{
\begin{array}{rl}
S^*_\zs{l}&=\frac{1}{H_\zs{l}}\,
\int^{\tau_l}_\zs{t_0}\,Q\left(\frac{y_s-x_\zs{l}}{h}\right)\,\d y_\zs{s}\,,\\[5mm]
\tau_\zs{l}&=\,\inf\{t\ge t_0\,:\,\int^t_\zs{t_0}\,
Q\left(\frac{y_s-x_\zs{l}}{h}\right)\,\d s
\ge\,H_\zs{l}\}\,,\ 
\end{array}
\right.
\end{equation}
where  $h=(b-a)/(2n)$,
 $Q(z)=\Chi_\zs{\{|z|\le 1\}}$ and
$$
H_l=(T-t_0)(2\tilde{q}_T(x_l)-\e^2_\zs{T})h
$$
with
$$
\tilde{q}_T(x_l)=\max\{\hat{q}(x_l)\,,\,\epsilon_\zs{T}\}
\quad\mbox{and}\quad
 \hat{q}(x_l)\,=\,
\frac{1}{2t_0h}\,
\int^{t_0}_\zs{0}\,Q\left(\frac{y_s-x_\zs{l}}{h}\right)\,\d s\,.
$$
Note that $\tau_\zs{l}<\infty$ a.s., for any $S\in \Sigma_\zs{L}$ and 
 for all $1\le l\le n$ (see, \cite{GaPe3}).

Moreover, we assume that the parameters
$t_0=t_0(T)$ and $\epsilon_\zs{T}$ satisfy the following conditions

\begin{itemize}
\item[$\H_\zs{1})$]
{\sl For any $T\ge 32$,
$$
16\le t_0\le T/2
\quad\mbox{and}\quad
\sqrt{2}/t^{1/8}_0\le\epsilon_\zs{T}\le 1\,.
$$
}
\item[$\H_\zs{2})$]
{\sl For any $\delta>0$ and $\nu>0$,
$$
\lim_\zs{T\to\infty}\,
T^{\nu}\,e^{-\delta\sqrt{t_0}}\,=0\,.
$$
}
\item[$\H_\zs{3})$]
$$
\lim_\zs{T\to\infty}\,t_0(T)\,=\,\infty\,,\quad
\lim_\zs{T\to\infty}\,\e_\zs{T}\,=\,0\,,\quad
\lim_\zs{T\to\infty}\,T\e_\zs{T}/t_0(T)\,=\,\infty\,.
$$
\end{itemize}
>From \eqref{1.1}\,,\eqref{AD.1}--\eqref{AD.1-2} we obtain the discrete regression model
$$
S^*_\zs{l}\,=\,S(x_l)\,+\,\zeta_l\,.
$$
The error term $\zeta_l$ is represented as the following sum of the approximated term
$B_l$ and the stochastic term
$$
\zeta_l\,=\,B_l\,+\,\frac{1}{\sqrt{H_\zs{l}}}\,\xi_\zs{l}\,,
$$
where 
\begin{align*}
B_\zs{l}&=\,\frac{1}{H_l}\,
\int^{\tau_l}_\zs{t_0}\,Q\left(\frac{y_s-x_\zs{l}}{h}\right)\,(S(y_s)\,-\,S(x_l))\d s\,,\\
\xi_\zs{l}&=\,\frac{1}{\sqrt{H_l}}\,
\int^{\tau_l}_\zs{t_0}\,Q\left(\frac{y_s-x_\zs{l}}{h}\right)\,\d w_s\,.
\end{align*}
Moreover, note that for any function $S\in W^{k}_\zs{r}$ 
\begin{equation}\label{AD.2}
|B_l|\,\le\,2L\,h\,=L\,(b-a)/n\,.
\end{equation}
It is easy to see that the random variables
$(\xi_l)_\zs{1\le l\le n}$ are i.i.d. normal $\cN(0,1)$. 
Therefore, by putting  
$$
\Gamma=\Gamma_\zs{T}=\{\max_\zs{1\le l\le n}\tau_l\,\le\,T\}
\quad\mbox{and}\quad
Y_l\,=\,S^*_\zs{l}\,\Chi_\zs{\Gamma}\,,
$$
we obtain on the set $\Gamma$ the following regression model 
\begin{equation}\label{AD.3}
Y_l\,=\,S(x_\zs{l})\,+\,\zeta_l\,, \ \ \ \zeta_l\,=\,B_\zs{l}\,+\,\sigma_l\,\xi_l\,,
\end{equation}
where
$$
\sigma^2_l\,=\,\frac{n}{(T-t_0)(\tilde{q}_T(x_l)-\e^2_\zs{T}/2)(b-a)}
\le \frac{4}{\e_\zs{T}(b-a)}:=\sigma_\zs{*}\,.
$$
In Appendix A.1 we prove  
the following result:
\begin{proposition}\label{Pr.AD.1}
Suppose that the parameters $t_0$ and $\epsilon_\zs{T}$ satisfy the conditions 
$\H_\zs{1})$--$\H_\zs{3})$. Then, for any $L\ge 1$,
% and $R\ge 1$,
$$
\lim_\zs{T\to\infty}\,
\sup_\zs{S\in\Sigma_\zs{L}}\,\sup_\zs{1\le l\le n}
\,\E_\zs{S}\,|\ov{\sigma}_\zs{l}|\,=\,0\,,
$$
where 
$$
\ov{\sigma}_\zs{l}= \sigma^2_\zs{l}\,-\,\frac{1}{q_\zs{S}(x_l)(b-a)}\,.
$$
\end{proposition}
Now we suppose that the parameters $k$ and $r$ of the space $W_r^k$ in \eqref{2.6}
are unknown. We describe the adaptive procedure from \cite{GaPe6}.
First we fixe
$\varepsilon>0$ and we define the sieve $\cA_\zs{\ve}$ in the space
$\bbn\times \bbr_\zs{+}$ :
\begin{equation}\label{AD.4}
\cA_\zs{\ve}=\{1,\ldots,k_*\}\times\{t_1,\ldots,t_m\}\,,
\end{equation}
where $k_*=[1/\sqrt{\ve}]$, $t_i=i\ve$, $m=[1/\ve^2]$ and we take $\ve=1/\ln n$.
 We remind of that $n=2[(T-1)/2]+1\ge 30$, due to condition $\H_\zs{1})$.

For any $\alpha=(\beta,t)\in\cA_\zs{\ve}$ we define the weight vector 
$\lambda_\zs{\alpha}=(\lambda_\zs{\alpha}(1),\ldots,\lambda_\zs{\alpha}(n))^{\prime}$ with
\begin{equation}\label{AD.5}
\lambda_\zs{\alpha}(j)\,=\,\left\{
\begin{array}{cc}
1\,,&\quad  \mbox{for}\quad 1\le j\le j_0\,,\\[2mm]
\left(1-(j/\omega_\alpha)^\beta\right)_\zs{+}&\,,\quad \mbox{for}\quad j_0<j\le n\,, 
\end{array}
\right.
\end{equation}
where $j_0=j_0(\alpha)=\left[\omega_\zs{\alpha}/\ln(n+2)\right]+1$, 
$$
\omega_\zs{\alpha}=(A_\zs{\beta}\,t\,n)^{1/(2\beta+1)}
\quad\mbox{and}\quad
A_\zs{\beta}=\frac{(b-a)^{2\beta+1}(\beta+1)(2\beta+1)}{\pi^{2\beta}\beta}\,.
$$
For any $\alpha\in \cA_\zs{\ve}$, through the weight 
$\lambda_\zs{\alpha}=(\lambda_\zs{\alpha}(1),\ldots,\lambda_\zs{\alpha}(n))'$ we construct
 the weighted least squares estimator

\begin{equation}\label{AD.6}
\left\{
\begin{array}{cl}
\hat{S}_\zs{\alpha}&\,=\,S_\zs{0}+\sum^n_\zs{j=1}\,\lambda_\zs{\alpha}(j)\,\hat{\theta}_\zs{j,n}\,\phi_\zs{j}\,
\Chi_\zs{\Gamma}\,,\\[5mm]
\hat{\theta}_\zs{j,n}&\,=\,
((b-a)/n)\,\sum^n_\zs{l=1}\,(Y_l-S_\zs{0}(x_\zs{l}))\,\phi_j(x_l)\,.
\end{array}
\right.
\end{equation}
\vspace{4mm}

We remind of (see Section 4 in \cite{GaPe6}) that to construct an adaptive procedure one has to minimize
the empiric squared error of estimator \eqref{AD.6} over the weight family 
$\{\lambda_\zs{\alpha}\,,\,\alpha\in\cA_\zs{\ve}\}$. A difficulty appears since the empiric squared error 
contains a term which depends on unknown function $S$. We estimate this term as follows
$$
\tilde{\theta}_\zs{j,n}=
\hat{\theta}^2_\zs{j,n}-\,\frac{(b-a)^2}{n}\,s_\zs{j,n}
\quad\mbox{with}\quad
s_\zs{j,n}\,=\,\frac{1}{n}\,\sum^n_\zs{l=1}\,\sigma^2_l\,\phi^2_\zs{j}(x_l)\,.
$$
For any $\lambda\in \{\lambda_\zs{\alpha}\,,\alpha\in\cA_\zs{\ve}\}$
we define the empiric cost function $J_\zs{n}(\lambda)$ by the following way
$$
J_\zs{n}(\lambda)\,=\,\sum^n_\zs{j=1}\,\lambda^2(j)\hat{\theta}^2_\zs{j,n}\,-
2\,\sum^n_\zs{j=1}\,\lambda(j)\,\tilde{\theta}_\zs{j,n}\,
+\,\frac{1}{\ln T} P_\zs{n}(\lambda)
$$
with the penalty term defined as
$$
P_\zs{n}(\lambda)=\frac{|\lambda|^2(b-a)^2s_\zs{n}}{n}\,,
$$
where $|\lambda|^2=\sum^n_\zs{j=1} \lambda^2(j)$ and $s_\zs{n}=n^{-1}\,\sum^n_\zs{l=1}\,\sigma^2_l$.
We set
\begin{equation}\label{AD.8}
\hat{\alpha}=\mbox{agrmin}_\zs{\alpha\in\cA_\zs{\ve}}\,J_n(\lambda_\zs{\alpha})
\quad\mbox{and}\quad
\hat{S}_\zs{*}=\hat{S}_\zs{\hat{\alpha}}\,.
\end{equation}
In \cite{GaPe6} we proved the following  non-asymptotic oracle inequality.
\begin{theorem}\label{Th.3.1}
Assume that $S\in \Sigma_\zs{L}$ and $\hat{S}_\zs{\alpha}$ is defined in  \eqref{AD.6}.
Then, for any $T\ge 32$, 
the adaptive estimator \eqref{AD.8}
 satisfies the following inequality 
\begin{equation}\label{AD.9}
\cR(\hat{S}_\zs{*},S)
\,\le\,(1+D(\rho))
\,\min_\zs{\alpha\in\cA_\zs{\ve}}\,\cR(\hat{S}_\zs{\alpha},S)\,+\,
\frac{\cB_\zs{T}(\rho)}{n}\,,
\end{equation}
where 
$$
\rho=1/(6+\ln n)
\quad\mbox{and}\quad
n=2[(T-1)/2]+1\,.
$$
 Moreover,
 the functions $D(\rho)$ and $\cB_\zs{T}(\rho)$ defined in 
Theorem 4.2 from \cite{GaPe6} are such that
$\lim_\zs{\rho\to 0}D(\rho)=0$ 
and, for any $\delta>0$ ,
\begin{equation}\label{AD.10}
\lim_\zs{T\to \infty}
\frac{\cB_\zs{T}(\rho)}{T^\delta}=0\,.
\end{equation}
\end{theorem}
Our principal goal in this paper is to show that the inequality \eqref{AD.9}
is sharp in asymptotic sense, i.e. it yields inequalities \eqref{2.9} and
\eqref{2.10}.

\medskip
\section{Upper bound}\label{Se.Up}
\subsection{Known smoothness}

 We start with the estimation problem (\ref{1.1})
under the condition that $S\in W^{k}_r$ with known
 parameters $k$, $r$ and $J(S)$ defined in \eqref{2.7}.
In this case we use the estimator from family \eqref{AD.6}
\begin{equation}\label{U.1}
\tilde{S}=\hat{S}_\zs{\tilde{\alpha}}
\quad\mbox{with}\quad
\tilde{\alpha}=(k,\tilde{t}_n)\,,\quad \tilde{t}_\zs{n}=
%t_\zs{\tilde{l}}=
\tilde{l}_\zs{n}\ve\,,
\end{equation}
 where 
$$
\tilde{l}_\zs{n}=\inf\{i\ge 1\,:\,i\ve\ge \ov{r}(S)\}\,,
%\quad\mbox{and}
\quad \ov{r}(S)=r/J(S)
$$
and $\ve=\ve_n=1/\ln n$.
Note that for sufficiently large $T$, therefore large 
$m=[1/\ve^2]=[\ln^2 n]$, the parameter 
$\tilde{\alpha}$ belongs to the set \eqref{AD.4}.
In this section we obtain the upper bound for
the empiric squared error of the estimator \eqref{U.1}. We define the empiric squared error 
of the estimator $\tilde{S}$ as
$$
\|\tilde{S}-S\|^2_\zs{n}=\frac{b-a}{n}\sum^n_\zs{l=1}(\tilde{S}(x_\zs{l})-S(x_\zs{l}))^2\,,
$$
where the points $(x_\zs{l})_\zs{1\le l\le n}$ are defined in \eqref{AD.1}.
\begin{theorem}\label{Th.Up.1} 
The estimator $\tilde{S}$ satisfies the following asymptotic upper bound
\begin{equation}\label{U.2}
\limsup_{T\to\infty}\,T^{2k/(2k+1)}\,
\sup_{S\in W^{k}_r}
\frac{1}{\gamma(S)}\,\E_\zs{S}\|\tilde{S}-S\|^2_\zs{n}\,\Chi_\zs{\Gamma}\,\le 1\,.
\end{equation}
\end{theorem}
\noindent {\bf Proof.} We denote $\tilde{\lambda}=\lambda_\zs{\tilde{\alpha}}$
and $\tilde{\omega}=\omega_\zs{\tilde{\alpha}}$.
Now we remind of  that the sieve Fourier coefficients  $(\hat{\theta}_\zs{j,n})$
defined in \eqref{AD.6} 
 satisfy on the set $\Gamma$
the following relation (see \cite{GaPe6})
\begin{equation}\label{U.3}
\hat{\theta}_\zs{j,n}\,=\,\theta_\zs{j,n}\,+\,\zeta_\zs{j,n}
\end{equation}
with
$$
\theta_\zs{j,n}=
\frac{b-a}{n}\,\sum^n_\zs{l=1}\,(S(x_\zs{l})-S_\zs{0}(x_\zs{l}))\,\phi_j(x_l)
$$
and
$$
\zeta_\zs{j,n}=(\zeta,\phi_j)_\zs{n}\,=\,\frac{b-a}{\sqrt{n}}\,\xi_\zs{j,n}\,+\,\delta_\zs{j,n}\,,
$$
where
\begin{equation}\label{U.5}
\xi_\zs{j,n}\,=\,\frac{1}{\sqrt{n}}\,\sum^n_\zs{l=1}\,\sigma_l\,\xi_l\,\phi_j(x_l)\ \ \ \mbox{and}\ \ \ \
\delta_\zs{j,n}\,=
\,\frac{b-a}{n}\,\sum^n_\zs{l=1}\,B_l\,\phi_j(x_l)\,.
\end{equation}
The inequality \eqref{AD.2} implies that
\begin{equation}\label{U.6}
|\delta_\zs{j,n}|\le L\,(b-a)^{3/2}/n\,.
\end{equation}
On the set $\Gamma$
we can represent the empiric squared error as follows 
\begin{align}\nonumber
\|\tilde{S}-S\|_n^2
&=\sum_{j=1}^{n}\,(1\,-\,\tilde{\lambda}(j))^2\,\theta^2_\zs{j,n}\,+2(b-a)M_\zs{n}\\ \label{U.7}
&+
2\sum_{j=1}^{n}\,(1\,-\,\tilde{\lambda}(j))\,\tilde{\lambda}(j)\,\theta_{j,n}\,\delta_\zs{j,n}
+\,
\sum_{j=1}^{n}\,\tilde{\lambda}^2(j)\,\zeta^2_\zs{j,n}\,,
\end{align}
where
$$
M_\zs{n}\,=\,\frac{1}{\sqrt{n}}
\sum_{j=1}^{n}\,(1\,-\,\tilde{\lambda}(j))\,\tilde{\lambda}(j)\,\theta_{j,n}\,\xi_\zs{j,n}\,.
$$
Note that, for any $\rho>0$,
\begin{align*}
2\sum_{j=1}^{n}\,(1\,-\,\tilde{\lambda}(j))\,\tilde{\lambda}(j)\,\theta_{j,n}\,\delta_\zs{j,n}
\le 
\rho\,\sum_{j=1}^{n}(1-\tilde{\lambda}(j))^2\,\theta^2_{j,n}\,
+
\rho^{-1}\sum_{j=1}^{n}\tilde{\lambda}^2(j)\,\delta^2_\zs{j,n}\,.
\end{align*}
Therefore by \eqref{U.5}-\eqref{U.6} we obtain that
\begin{align*}
\|\tilde{S}-S\|_n^2
&\le (1+\rho)\sum_{j=1}^{n}\,(1\,-\,\tilde{\lambda}(j))^2\,\theta^2_\zs{j,n}\,+\,2(b-a)M_\zs{n}\\ 
&
+\frac{L^2(b-a)^3}{\rho n}+
\sum_{j=1}^{n}\,\tilde{\lambda}^2(j)\,\zeta^2_\zs{j,n}\,.
\end{align*}
By the same way we estimate the last term in the right-hand part as
\begin{align*}
\sum_{j=1}^{n}\,\tilde{\lambda}^2(j)\,\zeta^2_\zs{j,n}&\le\,
\frac{(1+\rho)(b-a)^2}{n}\sum_{j=1}^{n}\,\tilde{\lambda}^2(j)\,\xi^2_\zs{j,n}\\
&+\,
(1+\rho^{-1})\,\frac{L^2(b-a)^3}{n}\,.
\end{align*}
Therefore on the set $\Gamma$ we find that
\begin{align}\nonumber
\|\tilde{S}_n-S\|_n^2&\le(1+\rho)\hat{\gamma}_\zs{n}(S)+2(b-a)M_\zs{n}\\ \label{U.8}
&+(1+\rho)\Delta_\zs{n}
+\frac{L^2(b-a)^3(\rho+2)}{\rho n}\,,
\end{align}
where
\begin{align*}
\hat{\gamma}_\zs{n}(S)\,&=\,\sum_{j=1}^{n}\,(1\,-\,\tilde{\lambda}(j))^2\,\theta^2_\zs{j,n}\,+\,
\frac{J(S)}{(b-a)n}\sum_{j=1}^{n}\,\tilde{\lambda}^2(j)\,,\\
%\frac{J(S)}{(b-a)n}\sum_{j=1}^{n}\,\tilde{\lambda}^2(j)\,,\\
\Delta_\zs{n}\,&=\,
\frac{1}{n}\sum_{j=1}^{n}\,\tilde{\lambda}^2(j)\,\left((b-a)^2\xi^2_\zs{j,n}-\frac{J(S)}{b-a}\right)\,.
\end{align*}
Let us estimate the first term in the right-hand part of \eqref{U.8}.
Note that the bounds \eqref{2.4} imply the corresponding bounds for the function 
$J(S)$, i.e.
\begin{equation}\label{U.9}
0<\frac{b-a}{ q^*}\le\inf_\zs{S\in\Sigma_\zs{L}}\,J(S)
\le
\sup_\zs{S\in\Sigma_\zs{L}}\,J(S)\le \frac{b-a}{q_\zs{*}}<\infty\,.
\end{equation}
This implies directly that
\begin{equation}\label{U.10}
\lim_\zs{n\to\infty} \sup_\zs{S\in\Sigma_\zs{L}}\,\left|\frac{\tilde{t}_\zs{n}}{\ov{r}(S)}-1\right|=0\,.
\end{equation}
Moreover, from the definition \eqref{AD.5} we get that
\begin{equation}\label{U.11}
\lim_{n\to\infty}\sup_\zs{S\in\Sigma_\zs{L}}\left|
\frac{\sum_{j=1}^n\,\tilde{\lambda}(j)^2}{n^{1/(2k+1)}}
-
\frac{(A_\zs{k}\ov{r}(S))^{1/(2k+1)}\,
k^2}{(k+1)(2k+1)}
\right|=0\,.
\end{equation}
Taking this into account, in Appendix~\ref{Su.A.2} we show that
\begin{equation}\label{U.12}
\limsup_\zs{T\to\infty}\sup_\zs{S\in W^{k}_r}\,T^{2k/(2k+1)}
\frac{\hat{\gamma}_\zs{n}(S)}{\gamma(S)}\,\le\,1\,.
\end{equation}
To estimate the second term in the right-hand part of inequality \eqref{U.8}
we use Lemma~\ref{Le.A.1}. We get
$$
\E_\zs{S}\,M^2_n\le\,\frac{\sigma_*}{(b-a)n}\,
\sum^{n}_\zs{j=1}\,\theta^2_\zs{j,n}=\frac{\sigma_*}{(b-a) n}
\|S\|^2_\zs{n}
\le \frac{\sigma_*\nu^*}{(b-a) n}
\,,
$$
where the constants $\sigma_*$ and $\nu^*$ are defined in \eqref{AD.3} and
\eqref{2.2-1}, respectively.

Taking into account that $\E_\zs{S}\,M_n=0$
and  making use of Proposition 3.1 from \cite{GaPe6}
we obtain that
\begin{align*}
|\E_\zs{S}\,M_n\,\Chi_\zs{\Gamma}|\,=\,|\E_\zs{S}\,M_n\,\Chi_\zs{\Gamma^c}|\,
\le\sqrt{\sigma_*\Pi_\zs{T}}\,\frac{L}{\sqrt{n}}\,.
%\le\sqrt{\sigma_*\Pi_\zs{T}}\,\frac{L}{(b-a)\sqrt{n}}\,.
\end{align*}
Therefore
\begin{equation}\label{U.13}
\lim_\zs{T\to\infty}T^{2k/(2k+1)}\sup_\zs{S\in W^{k}_r}|\E_\zs{S}\,M_n\,\Chi_\zs{\Gamma}|=0\,.
\end{equation}
Now we show that
\begin{equation}\label{U.14}
\lim_\zs{T\to\infty}T^{2k/(2k+1)}\sup_\zs{S\in W^{k}_r}|\E_\zs{S}\,\Delta_n|=0\,.
\end{equation}
First of all, note that, for $j\ge 2$,
\begin{align}\nonumber
(b-a)^2\E_\zs{S}\,\xi^2_\zs{j,n}&=\frac{(b-a)^2}{n}\E_\zs{S}\sum^n_\zs{l=1}
\sigma^2_\zs{l}\phi_j^2(x_\zs{l})\\ \label{U.14-1}
&=(b-a)\E_\zs{S}\,s_\zs{n}+(b-a)\E_\zs{S}\ov{\varsigma}_\zs{j,n}\,,
\end{align}
where 
$$
\ov{\varsigma}_\zs{j,n}=
\frac{1}{n}
\sum^n_\zs{l=1}\sigma^2_\zs{l}\ov{\phi}_j(x_\zs{l})
\quad\mbox{with}\quad
\ov{\phi}_j(z)=(b-a)\phi^2_\zs{j}(z)-1\,.
$$
Moreover, 
\begin{align*}
\sup_\zs{S\in W^k_\zs{r}}\,\left|(b-a)\E_\zs{S}\,s_n-\frac{J(S)}{b-a}\right|\,
&\le\,\frac{b-a}{n}\,\sum^n_\zs{l=1}
\sup_\zs{S\in W^k_\zs{r}}\,\E_\zs{S}\,|\ov{\sigma}_\zs{l} |\\
&+\,
\sup_\zs{S\in W^k_\zs{r}}\,
\left|\frac{1}{b-a}\int^b_a\,q^{-1}_\zs{S}(x)\,\d x-\frac{1}{n}\sum^n_\zs{l=1}\, q^{-1}_\zs{S}(x_l)\right|\\
&\le (b-a)
\sup_\zs{S\in W^k_\zs{r}}\max_\zs{1\le l\le n}\,\E_\zs{S}|\ov{\sigma}_l|+
\frac{q'_\zs{*}(b-a)}{q_\zs{*}n}\,,
\end{align*}
where $q'_\zs{*}=\max_\zs{a\le x\le b}\,\sup_\zs{S\in\Sigma_\zs{L}}\,|q'_\zs{S}(x)|$.
Therefore
Proposition~\ref{Pr.AD.1} implies that
$$
\lim_\zs{n\to\infty}\,T^{2k/(2k+1)}
\sup_\zs{S\in W^k_\zs{r}}\,\left|(b-a)\E_\zs{S}\,s_n-\frac{J(S)}{b-a}\right|=0.
$$
To estimate the second term in \eqref{U.14-1} 
we make use of  Lemma 6.2 from \cite{GaPe4}. We have
\begin{align*}
\left|\sum_{j=1}^{n}\,\tilde{\lambda}^2(j)\,\ov{\varsigma}_\zs{j,n}\right|\,&=\,
\frac{1}{n}\left|\sum_{l=1}^{n}\,\sigma^2_l\,
\sum_{j=1}^{n}\,\tilde{\lambda}^2(j)\,\ov{\phi}_\zs{j}(x_l)\right|\\
&\le\,\frac{1}{n}\,\sum_{l=1}^{n}\,\sigma^2_l\,
\left|\sum_{j=1}^{n}\,\tilde{\lambda}^2(j)\,\ov{\phi}_\zs{j}(x_l)\right|\\
&\le\,
\sigma_*\,(2^{2k+1}+2^{k+2}+1)
\le\,
5\sigma_*\,2^{2k}
\quad\quad\mbox{a.s..}
\end{align*}
Thus from \eqref{U.11} we obtain \eqref{U.14}.
Moreover, we can calculate that
$$
\E_\zs{S}\,\xi^4_\zs{j,n}\le \frac{3\sigma^2_\zs{*}}{(b-a)^2}\,.
$$
Due to Proposition 3.1 from \cite{GaPe6}, we obtain that
\begin{align*}
\E_\zs{S}|\Delta_\zs{n}|\Chi_\zs{\Gamma^c}&\le 
\frac{(b-a)^2}{n}\sum^n_\zs{j=1}\E_\zs{S}\xi^2_\zs{j,n}\Chi_\zs{\Gamma^c}+
J(S)\Pi_\zs{T}\\
&\le \frac{\sqrt{3}\sigma_\zs{*}}{b-a}\sqrt{\Pi_\zs{T}}+\frac{1}{q_\zs{*}}\Pi_\zs{T}\,.
\end{align*}
This means that
$$
\lim_\zs{n\to\infty}\,T^{2k/(2k+1)}
\sup_\zs{S\in W^k_\zs{r}} \E_\zs{S}|\Delta_\zs{n}|\Chi_\zs{\Gamma^c}
=0\,.
$$
Therefore by \eqref{U.14} we get that
$$
\lim_\zs{n\to\infty}\,T^{2k/(2k+1)}
\sup_\zs{S\in W^k_\zs{r}} |\E_\zs{S}\,\Chi_\zs{\Gamma}\Delta_\zs{n}|
=0\,.
$$
Hence Theorem~\ref{Th.Up.1}.
\endproof

\subsection{Unknown smoothness}
In this subsection we prove Theorem~\ref{Th.2.1}.
First of all notice that inequalities \eqref{U.9} yield
$$
\inf_{S\in W^{k}_r}\,\gamma(S)\,>\,0\,.
$$
Therefore
 Theorem~\ref{Th.Up.1}, upper bound \eqref{2.2-1} and Proposition 3.1 from \cite{GaPe6} imply that
\begin{equation}\label{P.1}
\limsup_{T\to\infty}\,T^{2k/(2k+1)}\,
\sup_{S\in W^{k}_r}
\frac{1}{\gamma(S)}\,\E_\zs{S}\|\tilde{S}-S\|^2_\zs{n}\,\le 1\,.
\end{equation}
Let us remind of of that we define the estimator  $\tilde{S}$ from the sieve \eqref{AD.1} onto all
interval $[a,b]$ by the standard method as
\begin{equation}\label{P.1-1}
\tilde{S}(x)=\tilde{S}(x_\zs{1})\Chi_\zs{\{a\le x\le x_\zs{1}\}}+
\sum^n_\zs{l=2}\tilde{S}(x_\zs{l})\Chi_\zs{\{x_\zs{l-1}< x\le x_\zs{l}\}}\,,
\end{equation}
where $\Chi_A$ is the indicator of a set $A$.
Putting $\varrho(x)=\tilde{S}(x)-S(x)$ we find that 
\begin{align*}
\|\varrho\|^2&=\|\varrho\|^2_\zs{n}+
2\sum^n_\zs{l=1}\int^{x_\zs{l}}_\zs{x_\zs{l-1}}\varrho(x_\zs{l})(S(x_\zs{l})-S(x))\d x\\
&+\sum^n_\zs{l=1}\int^{x_\zs{l}}_\zs{x_\zs{l-1}}(S(x_\zs{l})-S(x))^2\d x\,.
\end{align*}
For any $0<\epsilon<1$,
we estimate the norm $\|\varrho\|^2$ as
\begin{align*}
\|\varrho\|^2\le(1+\epsilon)\|\varrho\|^2_\zs{n}+
(1+\epsilon^{-1})\sum^n_\zs{l=1}\int^{x_\zs{l-1}}_\zs{x_\zs{l-1}}(S(x_\zs{l})-S(x))^2\d x\,.
\end{align*}
This means that, for any $S\in \Sigma_\zs{L}$,
\begin{equation}\label{P.2}
\cR(\tilde{S},S)\,\le(1+\epsilon)\E_\zs{S}\|\varrho\|^2_\zs{n}+
(1+\epsilon^{-1})\,\frac{L^2(b-a)^3}{n^2}\,.
\end{equation}
We recall that $\tilde{S}=\hat{S}_\zs{\tilde{\alpha}}$ with $\tilde{\alpha}\in \cA_\zs{\ve}$.
Therefore, Theorem~\ref{Th.3.1} with inequalities \eqref{P.1}--\eqref{P.2} imply Theorem~\ref{Th.2.1}
\endproof
\medskip

\section{Lower bound}\label{Se.Lo}

In this section we prove Theorem~\ref{Th.2.2}. 
We folow the proof of Theorem 4.2 from 
\cite{GaPe8}. Similarly, we start with the approximation for an indicator
function, i.e. for any
for $\eta>0$,  we set
\begin{equation}\label{Lo.1}
I_\zs{\eta}(x)
=\eta^{-1}\int_{\bbr}\Chi_\zs{(|u|\le 1-\eta)}G\left(\frac{u-x}{\eta}\right)\,\d u\,,
\end{equation}
where the kernel $V\in\C^{\infty}(\bbr)$  is 
a probability density on $[-1,1]$.
It is easy to see that  $I_\zs{\eta}\in \cC^{\infty}$ and for any $m\ge 1$ and
any integrable function $f(x)$
$$
\lim_\zs{\eta\to 0}\,
\int_{\bbr}f(x)I^m_{\eta}(x)\,\d x\,=\,\int_{-1}^{1}f(x)\,\d x\,.
$$
Further, we will make use of the following trigonometric basis 
$\{e_j, j\ge 1\}$ in $\cL_2[-1,1]$ with
\begin{equation}\label{Lo.2}
e_1(x)=1/\sqrt{2}\,,\  e_j(x)\,=\,Tr_j(\pi[j/2] x)\,,\ j\ge 2\,.
\end{equation}
Here $Tr_\zs{l}(x)=\cos(x)$ for even $l$ and 
$Tr_\zs{l}(x)=\sin(x)$ for odd $l$.

Moreover, we denote 
$$
J_0=J(S_\zs{0})\,,\quad q_0=q_\zs{S_\zs{0}}\,,\quad \gamma_0=\gamma(S_\zs{0})\,,
$$
where the function $S_\zs{0}$ is defined in \eqref{2.5}.

Let us now fixe some arbitrary $0<\ve<1$ and according to \cite{GaPe8}
we put
\begin{equation}\label{Lo.3}
h=(\upsilon^{*}_\zs{\ve})^{\frac{1}{2k+1}}\,N_\zs{T}\,T^{-\frac{1}{2k+1}}
\end{equation}
with 
$$
\upsilon^{*}_\zs{\ve}=\frac{k\pi^{2k}J_\zs{0}}{(1-\ve)r2^{2k+1}(k+1)(2k+1)}
\quad\mbox{and}\quad
N_\zs{T}=\ln^4 T\,.
$$

To construct a parametric family we divide the interval
$[a,b]$ by the intervals $[\tilde{x}_m-h\,,\tilde{x}_m+h]$ with
$\tilde{x}_m=a+2hm$. The maximal number of such intervals is equal to
$$
M=[(b-a)/(2h)]-1\,.
$$
Onto each interval $[\tilde{x}_m-h\,,\tilde{x}_m+h]$, we approximate
any unknown function by a trigonometric series with $N$ terms, i.e.
for any array $z=(z_\zs{m,j})_\zs{ 1\le m\le M\,, 1\le j\le N}$,
 we set
\begin{equation}\label{Lo.4}
S_\zs{z,T}(x)=S_\zs{0}(x)+\sum_{m=1}^M\sum_{j=1}^N\,z_\zs{m,j}D_\zs{m,j}(x)
\end{equation}
with $D_\zs{m,j}(x)=e_j(v_m(x))I_\zs{\eta}\,(v_m(x))$
and $v_m(x)=(x-\tilde{x}_m)/h$.

Now to obtain the bayesian risk we  choose a prior distribution on $\bbr^{MN}$
by making use of the random  array 
$\theta=(\theta_\zs{m,j})_\zs{1\le m\le M\,, 1\le j\le N}$ defined as
\begin{equation}\label{Lo.5}
\theta_\zs{m,j}\,=\,t_\zs{m,j}\,\zeta_\zs{m,j}\,,
\end{equation}
where $\zeta_\zs{m,j}$ are i.i.d. gaussian $\cN(0,1)$ random variables and
the coefficients 
$$
t_\zs{m,j}=\frac{\sqrt{y^*_\zs{j}}}{\sqrt{Thq_0(\tilde{x}_m)}}\,.
$$
We chose the sequence $(y^*_\zs{j})_\zs{1\le j\le N}$ by thre samle way as in (8.11) 
in \cite{GaPe8}, i.e.
$$
y^*_j\,
=\,\Omega_\zs{T}\,j^{-k}-1
\quad\mbox{with}\quad
\Omega_\zs{T}=\frac{R^*_\zs{T}+\sum^{N}j^{2k}}{\sum^{N}j^{k}}
\,,
$$
where
$$
R^*_\zs{T}=\frac{J_\zs{0}k}{\hat{J}_\zs{0}(k+1)(2k+1)}\,N^{2k+1}\,,
$$
and
\begin{equation}\label{Lo.6}
\hat{J}_\zs{0}=2h\sum^{M}_\zs{m=1}\frac{1}{q_\zs{0}(\tilde{x}_\zs{m})}\,.
\end{equation}
In the sequel we make use of the following set

\begin{equation}\label{Lo.7}
\Xi_\zs{T}=\{\max_\zs{1\le m\le M}\,\max_\zs{1\le j\le N} |\zeta_\zs{m,j}|\le \ln T\}
\,.
\end{equation}
Obviously, that for any $p>0$
\begin{equation}\label{Lo.8}
\lim_\zs{T\to\infty}\,T^{p}\,\P(\Xi^{c}_\zs{T})\,=0\,.
\end{equation}
Note that on the set $\Xi_\zs{T}$
the uniform norm 
$$
|S_\zs{\theta,T}-S_\zs{0}|_\zs{*}=\sup_\zs{a\le x\le b}
|S_\zs{\theta,T}(x)-S_\zs{0}(x)|
$$ 
is bounded
\begin{equation}\label{Lo.9}
|S_\zs{\theta,T}-S_\zs{0}|_\zs{*}\,\le 
\frac{\ln T }{\sqrt{q_\zs{*}Th}}\sum^{N}_\zs{j=1}\,\sqrt{y^*_\zs{j}}
:=\epsilon_\zs{T}\,.
\end{equation}
Taking into account here that
\begin{equation}\label{Lo.10}
\lim_\zs{T\to\infty}\,\hat{J}_\zs{0}=J_\zs{0}
\end{equation}
it is easy to deduce that $\epsilon_\zs{T}\to 0$ as $T\to\infty$.

 For any estimator $\hat{S}_\zs{T}$, we denote by $\hat{S}^0_\zs{T}$ its projection on 
 $W^k_r$, i.e.\\
 $\hat{S}^0_\zs{T}=\hbox{\rm Pr}_\zs{W^k_r}(\hat{S}_\zs{T})$.
Since $W^k_\zs{r}$ is a convex set, we get that
$$
\|\hat{S}_\zs{T}-S\|^2\ge\|\hat{S}^0_\zs{T}-S\|^2\,.
$$
Therefore, denoting by $\mu_\zs{\theta}$  
 the distribution of $\theta$ in $\bbr^{d}$ with $d=MN$
 and
 taking into account \eqref{Lo.9}
we can write that
$$
\sup_\zs{S\in W^k_r}\,\frac{\cR(\hat{S}_\zs{T},S)}{\gamma(S)}
\ge\frac{1}{\gamma^*_\zs{T}}\,
\int_{\{z\in\bbr^d\,:\,S_\zs{z,n}\in W^k_\zs{r}\}\cap\Xi_\zs{T}}\,
\E_\zs{S_\zs{z,T}}\|\hat{S}^0_\zs{T}-S_\zs{z,T}\|^2\,\mu_{\vartheta}(\d z)
$$
with 
$$
\gamma^*_\zs{T}=\sup_\zs{S\in \cU_\zs{T}}\,\gamma(S)\,,
$$
where 
$$
\cU_\zs{T}=\{S\,:\,|S-S_\zs{0}|_\zs{*}\le \epsilon_\zs{T}\,,S(x)=S_\zs{0}(x)
\quad\mbox{for}\quad x\notin [a,b]\}\,.
$$
Since function \eqref{2.3} is continuous with respect to $S$, then
\begin{equation}\label{Lo.11}
\lim_\zs{T\to\infty}\gamma^*_\zs{T}=\gamma_\zs{0}\,.
\end{equation}

Making use of the distribution $\mu_\zs{\theta}$ we introduce
 the following Bayes risk
$$
\tilde{\cR}(\hat{S}_\zs{T})=
\int_\zs{\bbr^d}\cR(\hat{S}_\zs{T},S_\zs{z,T})
\mu_{\theta}(\d z)
$$
Now noting that $\|\hat{S}^0_\zs{T}\|^2\le r$ 
 through this risk we can write that
\begin{equation}\label{Lo.12}
\sup_\zs{S\in W_r^k}\,\frac{\cR(\hat{S}_\zs{T},S)}{\gamma(S)}\,
\ge\,\frac{1}{\gamma^*_\zs{T}}\,
\tilde{\cR}(\hat{S}^0_\zs{T})\,-\,\frac{2}{\gamma^*_\zs{T}}\,
\varpi_\zs{T}\,,
\end{equation}
with
$$
\varpi_\zs{T}=\E
(\Chi_\zs{\{S_\zs{\theta,T}\notin W^k_\zs{r}\}}\,+\,
\Chi_\zs{\Xi^c_\zs{T}})
(r+\|S_\zs{\theta,T}\|^2)\,.
$$
 Propostions 7.2--7.3 from \cite{GaPe8} imply that for any $p>0$
$$
\lim_\zs{T\to\infty}\,T^{p}\,\varpi_\zs{T}\,=0\,.
$$

Let us consider the first term in the right-hand side of \eqref{Lo.12}. 
To obtain a lower bound for this term we use the $\cL_2[a,b]$-orthonormal 
function family $(\tilde{e}_\zs{m,j})_\zs{1\le m\le M,1\le j\le N}$ which is defined as
$$
\tilde{e}_\zs{m,j}(x)=
\frac{1}{\sqrt{h}}e_j\left(v_m(x)\right)\Chi_\zs{\left(|v_m(x)|\le 1\right)}\,.
$$
We denote by $\hat{\lambda}_\zs{m,j}$ and $\lambda_\zs{m,j}(z)$ the Fourier coefficients
for the functions $\hat{S}^0_\zs{T}$ and $S_z$, respectively, i.e.
$$
\hat{\lambda}_\zs{m,j}=\int^b_a\,\hat{S}^0_\zs{T}(x)\tilde{e}_\zs{m,j}(x)\d x
\quad\mbox{and}\quad
\lambda_\zs{m,j}(z)=\int^b_a\,S_z(x)\tilde{e}_\zs{m,j}(x)\d x\,.
$$
Now it is easy to see that
$$
\|\hat{S}^0_\zs{T}-S_z\|^2 \ge \sum_{m=1}^M\sum_{j=1}^N\,
(\hat{\lambda}_\zs{m,j}-\lambda_\zs{m,j}(z))^2\,.
$$
Let us introduce the folowing $\cL_\zs{1}\to\bbr$ functional
$$
\ov{e}_\zs{j}(f)=
\int^1_{-1}\,e^2_\zs{j}(v)\,f(v)\,\d v\,
$$
Therefore from definition \eqref{Lo.4}
we obtain that
$$
\frac{\partial }{\partial z_\zs{m,j}}
\lambda_\zs{m,j}(z)
=\,\sqrt{h}\,\ov{e}_\zs{j}(I_\zs{\eta})\,.
$$

Now Lemma~\ref{Le.A.2} implies that
\begin{equation}\label{Lo.13}
\tilde{\cR}(\hat{S}^0_\zs{T})\,\ge\,h
\sum_{m=1}^M\sum_{j=1}^N\,\frac{\ov{e}^2_j(I_\zs{\eta})}
{(1+\varsigma_\zs{m,j}(T))\ov{e}_\zs{j}(I^2_\zs{\eta})\,q_\zs{0}(\tilde{x}_\zs{m})Th\,+\,t^{-2}_\zs{m,j}}\,,
\end{equation}
where
$$
\varsigma_\zs{m,j}(T)=
\E
\frac{\E_\zs{S_\zs{\theta,T}}\,\int^T_\zs{0}\,
D^2_\zs{m,j}(y_t)\,\d t
}
{Th \ov{e}_\zs{j}(I^2_\zs{\eta})\,q_\zs{0}(\tilde{x}_\zs{m})}
-1
\,.
$$
In Appendix we show that
\begin{equation}\label{Lo.14}
\lim_\zs{T\to\infty}
\max_\zs{1\le m\le M}\max_\zs{1\le j\le N}
\left|
\varsigma_\zs{m,j}(T)
\right|\,=\,0\,.
\end{equation}
Therefore taking this into account in inequality \eqref{Lo.13} we obtain that for sufficiently large 
$T$ and for arbitrary $\nu>0$ 
$$
\tilde{\cR}(\hat{S}^0_\zs{T})\,\ge\,
\frac{\hat{J}_\zs{0}}{2Th(1+\nu)}\,
\sum_{j=1}^N\,\tau_\zs{j}(\eta,y^*_\zs{j})\,,
$$
where
$$
\tau_\zs{j}(\eta,y)=\frac{\ov{e}^2_\zs{j}(I_\zs{\eta})y}{\ov{e}_\zs{j}(I^2_\zs{\eta})y+1}\,.
$$
By making use of limit equality (8.9) from \cite{GaPe8} we obtain that for sufficiently small
$\eta$ and sufficientlly large $T$
$$
\tilde{\cR}(\hat{S}^0_\zs{T})\,\ge\,\frac{1}{(1+\nu)^2}
\frac{\hat{J}_\zs{0}}{2Th}\,
\sum_{j=1}^N\,
\frac{y^*_\zs{j}}{y^*_\zs{j}+1}\,,
$$
where $\hat{J}_\zs{0}$ is defined in \eqref{Lo.7}. Thus making use of 
\eqref{Lo.10} this implies that
$$
\liminf_\zs{T\to\infty}\inf_\zs{\hat{S}_\zs{T}}\,T^{\frac{2k}{2k+1}}\,
\tilde{\cR}(\hat{S}_\zs{T})\,\ge\,(1-\ve)^{\frac{1}{2k+1}}\ \gamma_\zs{0}\,.
$$
Taking into account this inequelity in \eqref{Lo.12} and limit equality
\eqref{Lo.11} we obtain that for any $0<\ve<1$
$$
\liminf_\zs{T\to\infty}\inf_\zs{\hat{S}_\zs{T}}\,T^{\frac{2k}{2k+1}}\,
\sup_\zs{S\in W_r^k}\,\frac{\cR(\hat{S}_\zs{T},S)}{\gamma(S)}\,
\ge\,
(1-\ve)^{\frac{1}{2k+1}}\,.
$$
Taking here limit as $\ve\to 0$ implies Theorem~\ref{Th.2.2}.
\endproof

\medskip

\setcounter{section}{0}
\renewcommand{\thesection}{\Alph{section}}

\section{Appendix}\label{Se.A}
\subsection{Proof of Proposition~\ref{Pr.AD.1}}\label{Su.A.1}

We use all notations from \cite{GaPe6}. 
For any function $\psi\,:\,\bbr\to\bbr$ such that
\begin{equation}\label{A.1}
\sup_\zs{y\in\bbr}\,|\psi(y)|\,<\,\infty
\quad\mbox{and}\quad
\int^{+\infty}_\zs{-\infty}\,|\psi(y)|\,\d y\,\le\,c^*\,<\infty
\end{equation}
we set 
$$
\M_\zs{S}(\psi)\,=\,\int^{+\infty}_\zs{-\infty}\,\psi(y)\,q_\zs{S}(y)\,\d y
\quad\mbox{and}\quad
\Delta_\zs{T}(\psi)\,=\,\frac{1}{\sqrt{T}}\,\int^T_\zs{0}\,
(\psi(y_t)-\M_\zs{S}(\psi))\,\d t\,.
$$
In \cite{GaPe7} (see Theorem 3.2) we show that, for any $\nu>0$ and for any 
$\psi$ satisfying \eqref{A.1}, there exists $\gamma=\gamma(c^*,L)>0$ such that 
 the following inequality holds
\begin{equation}\label{A.2}
\sup_\zs{S\in\Sigma_\zs{L}}\,
\P_\zs{S}(|\Delta_\zs{T}(\psi)|\,\ge\,\nu)\,\le\,8\,e^{-\gamma \nu^2}\,.
\end{equation}
We shall apply this inequality to the function
$$
\psi_\zs{h,k}(y)=\frac{1}{h}\,Q\left(\frac{y-x_k}{h}\right)\,,
$$
for which $\int^{\infty}_\zs{-\infty}\psi_\zs{h,k}(y)\d y=2$.
Note now that
\begin{equation}\label{A.3}
2\hat{q}(x_k)\,-\,\M_\zs{S}(\psi_\zs{h,k})\,=\,\frac{1}{\sqrt{t_0}}\,
\Delta_\zs{t_0}(\psi_\zs{h,k})\,.
\end{equation}
Moreover,
$$
\M_\zs{S}(\psi_\zs{h,k})=\int^{1}_\zs{-1}q_\zs{S}(x_\zs{k}+hz)\d z\ge 2q_*\,,
$$
where $q_*$ is defined in \eqref{2.4}. Therefore we get that
\begin{align*}
\P_\zs{S}(\hat{q}(x_k)\,<\,\epsilon_\zs{T})\,&=\,
\P_\zs{S}\left(
\frac{1}{t_0}\,\int^{t_0}_0\,\psi_\zs{h,k}(y_t)\,\d t\,<\,2\epsilon_\zs{T}
\right)\\
&=\,
\P_\zs{S}\left(
\Delta_\zs{t_0}(\psi_\zs{h,k})\,<\,(2\epsilon_\zs{T}- \M_\zs{S}(\psi_\zs{h,k}))\sqrt{t_0}\right)\\
&\le\,
\P_\zs{S}\left(
\Delta_\zs{t_0}(\psi_\zs{h,k})\,<\,2\,(\epsilon_\zs{T}- q_*)\sqrt{t_0}\right)\,.
\end{align*}
Note that for $\epsilon_\zs{T}\le q_*/2$  the inequality
 \eqref{A.2} implies the following exponentielle upper bound
\begin{equation}\label{A.4}
\P_\zs{S}(\hat{q}(x_k)\,<\,
\epsilon_\zs{T})\,\le\,8\,e^{-\gamma q^2_*t_0}\,.
\end{equation}

Now we show that
\begin{equation}\label{A.5}
\lim_\zs{T\to\infty}\,
\sup_\zs{1\le l\le n}\,\sup_\zs{S\in\Sigma_\zs{L}}\,\frac{1}{\e_\zs{T}}\,
\E_\zs{S}\,\left|\tilde{q}_\zs{T}(x_l)\,-\,q_\zs{S}(x_\zs{l})\right|\,=\,0\,.
\end{equation}
To end this we have to prove that
\begin{equation}\label{A.6}
\lim_\zs{T\to\infty}\,
\sup_\zs{1\le l\le n}\,\sup_\zs{S\in\Sigma_\zs{L}}\,\frac{1}{\e_\zs{T}}\,
\E_\zs{S}\,\left|\hat{q}(x_l)\,-\,
\M_\zs{S}(\psi_\zs{h,l})/2
\right|\,=\,0\,.
\end{equation}
Indeed, from \eqref{A.2}--\eqref{A.3} we find
\begin{align*}
\E_\zs{S}\,\left|\hat{q}(x_l)\,-\,\M_\zs{S}(\psi_\zs{h,l})/2\right|\,&=\,
\frac{1}{\sqrt{t_0}}\,\E_\zs{S}\,|\Delta_\zs{t_0}(\psi_\zs{h,k})|\\
&=\,\frac{1}{\sqrt{t_0}}\,\int^\infty_0\,
\P_\zs{S}(|\Delta_\zs{t_0}(\psi_\zs{h,k})|\ge z)\,\d z\\
&\le\,\frac{8}{\sqrt{t_0}}\,
\int^\infty_\zs{0}\,e^{-\gamma z^2}\,\d z\,.
\end{align*}
The condition $\H_\zs{1})$ implies that $\e_\zs{T}\sqrt{t_\zs{0}}\to \infty$
as $T\to\infty$. Therefore this inequality implies \eqref{A.6}.
Moreover, taking into account that
 $h/\e_\zs{T}\to 0$ as $T\to\infty$ we obtain, for sufficiently large $T$,
the following bound
\begin{align*}
|\M_\zs{S}(\psi_\zs{h,l})/2-q_\zs{S}(x_l)|&\le\,
\int^1_\zs{-1}\,|q_\zs{S}(x_l+vh)\,-\,q_\zs{S}(x_l)|\,\d v\\
&\le\, q^{''}_*h^2\,\le \,\e^2_\zs{T}\,,
\end{align*}
where
$q^{''}_*=\sup_\zs{|x|\le R}\sup_\zs{S\in\Sigma_\zs{L}}\,|q^{''}_\zs{S}(x)|$.
>From this inequality, taking into account inequality \eqref{A.6} and the condition $\H_\zs{3})$,
we obtain \eqref{A.5}.

Since $T-2\le n\le T$, we find that,
for sufficiently large $T$ providing $\e_\zs{T}\le 1$, 
\begin{align*}
\E_\zs{S}\,\left|\sigma^2_l-\frac{1}{q_\zs{S}(x_l)(b-a)}\right|\,&=\,\frac{1}{b-a}
\E_\zs{S}\,\left|
\frac{2n}{(T-t_0)(2\tilde{q}(x_l)-\e^2_\zs{T})}\,-\,\frac{1}{q_\zs{S}(x_l)}
\right|\\
&\le\,
2\,\E_\zs{S}\,\frac{|\tilde{q}(x_l)\,-\,q_\zs{S}(x_l)|}{\e_\zs{T}\,q_\zs{*}(b-a)}\,+\,
\,\frac{\e_\zs{T}}{q_*(b-a)}\\
&+\,\frac{4}{(T-t_0)\e_\zs{T}(b-a)}
+\,\frac{2t_\zs{0}}{(T-t_0)\e_\zs{T}(b-a)}\,.
\end{align*}
The condition $\H_\zs{3})$ and \eqref{A.5} imply directly Proposition~\ref{Pr.AD.1}.
\endproof

\subsection{Proof of the limiting inequality \eqref{U.12}}\label{Su.A.2}

We set $\tilde{\iota}_\zs{0}=j_\zs{0}(\tilde{\alpha})$ and
$\tilde{\iota}_1=[\tilde{\omega} \ln(n+1)]$.
Then we can represent $\hat{\gamma}_\zs{n}(S)$
by the following way
$$
\hat{\gamma}_\zs{n}(S)=\sum_{j=\tilde{\iota}_\zs{0}}^{\tilde{\iota}_1}\,(1\,-\,\tilde{\lambda}(j))^2\,\theta^2_\zs{j,n}\,+\,
\frac{J(S)}{(b-a)n}\sum_{j=1}^{n}\,\tilde{\lambda}^2(j)+\Delta_\zs{1,n}
$$
with $\Delta_\zs{1,n}=\sum_{j=\tilde{\iota}_1}^{n}\,\theta^2_\zs{j,n}$.
Note now that,
for any  $0<\delta<1$,
\begin{align}\nonumber
\hat{\gamma}_\zs{n}(S)
&\,\le\,(1+\delta)
\sum_{j=\tilde{\iota}_0}^{\tilde{\iota}_1-1}\,(1\,-\,\tilde{\lambda}(j))^2\theta^2_\zs{j}\,+\,
\frac{J(S)}{(b-a)n}
\sum_{j=1}^n\,\tilde{\lambda}(j)^2\\ \label{A.7}
&+\,\Delta_\zs{1,n}+(1+1/\delta) \Delta_\zs{2,n}\,,
\end{align}
where $\Delta_\zs{2,n}=\,
\sum_{j=\tilde{\iota}_0}^{\tilde{\iota}_1-1}\,
(\theta_\zs{j,n}\,-\,\theta_\zs{j})^2$.

Due to the uniform convergence \eqref{U.10},
 Lemmas~6.1 and 6.3  from \cite{GaPe4} yield
$$
\lim_{n\to\infty}\sup_{S\in W^{k}_r}\,n^{2k/(2k+1)}\,\sum^2_{l=1}|\Delta_\zs{l,n}|=0\,.
$$
Now we set
$$
\upsilon_\zs{n}(S)=
n^{2k/(2k+1)}
\sup_\zs{j\ge \tilde{\iota}_\zs{0}}(1-\tilde{\lambda}(j))^2/\varpi_\zs{j}\,,
$$
with the sequence $\varpi_\zs{j}$  defined in \eqref{2.6} and
$$
\upsilon^*(S)=\left(\frac{b-a}{\pi}\right)^{2k}\,
\frac{1}{\left(A_\zs{k}\ov{r}(S)\right)^{2k/(2k+1)}}\,,
$$
where the coefficient $A_\zs{k}$ is defined in \eqref{AD.5}.
Moreover, one can calculate directly that
\begin{equation}\label{A.8}
\limsup_{T\to\infty}\sup_\zs{S\in\Sigma_\zs{L}}\frac{\upsilon_\zs{n}(S)}{\upsilon^*(S)}
\le\,1\,.
\end{equation}
Therefore, due to the definition \eqref{2.6} and to the fact that 
$$
\gamma(S)=\upsilon^*(S)r+\frac{J(S)}{b-a}
(A_\zs{k}\ov{r}(S))^{1/(2k+1)}\frac{2k^2}{(k+1)(2k+1)}\,,
$$
the inequality \eqref{A.7} and the limits \eqref{A.8} and \eqref{U.11} imply
 (\ref{U.12}).
\endproof

\medskip
\subsection{Moment bounds}\label{Su.A.3}

\begin{lemma}\label{Le.A.1}
Let $\xi_\zs{j,n}$ be defined in \eqref{U.5}. 
Then, for any real numbers $v_1,\ldots,v_n$,
\begin{align*}
\E\,\left(\sum_{j=1}^{n}\,v_j\,\xi_\zs{j,n}\right)^2\,\le\,
\frac{\sigma_*}{b-a}V_\zs{n}\,,
\quad
\E\,\left(\sum_{j=1}^{n}\,v_j\,\xi_\zs{j,n}\right)^4\,\le\,\frac{3\sigma^2_\zs{*}}{(b-a)^2}
V^2_\zs{n}\,,
\end{align*}
where $\sigma_*=\max_\zs{1\le j\le n}\,\sigma_j^2$ and $V_\zs{n}=\sum_{j=1}^{n}v_j^2$.
\end{lemma}
The proof of this Lemma is similar to the proof of Lemma 6.4 in \cite{GaPe4}.
\endproof

\medskip
\subsection{Application of the van Trees inequality to diffusion processes.}\label{Su.A.4}

Let 
$\left(\cC[0,T],\cB,(\cB_t)_\zs{0\le t\le T},(\P_\zs{\theta}\,,\theta\in\bbr^d)\right)$
 be a 
filtered statistical model with cylidric $\sigma$-fields 
$\cB_\zs{t}$
on $\cC[0,t]$ and $\cB=\cup_\zs{0\le t\le T}\cB_\zs{t}$.
As to the distributions $\P_\zs{\theta}$ we assume that 
it is  distribution in $\cC[0,T]$ of  the stochastic process $(y_t)_\zs{0\le t\le T}$  
governed by the stochastic differential equation
\begin{equation}\label{A.9}
\d y_t=S(y_t,\theta)\d t+\d w_t\,,\quad 0\le t\le T\,,
\end{equation}
where  
$\theta=(\theta_\zs{1},\ldots,\theta_\zs{d})'$ is vector of
unknown parameters, $w=(w_\zs{t})_\zs{0\le t\le T}$
is a standart Wiener process.
 Moreover, we assume also that $S$ is a linear function
with respect to $\theta$, i.e.
\begin{equation}\label{A.10}
S(y,\theta)=\sum^d_\zs{i=1}\,\theta_\zs{i}\,S_\zs{i}(y)\,,
\end{equation}
where
 the functions  $(S_\zs{i})_\zs{1\le i\le d}$ are bound and satisfy the Lipschitz
condition, i.e. for some constant $0<L<\infty$
$$
\max_\zs{1\le i\le d}\sup_\zs{x\in\bbr}\,|S_\zs{i}(x)|\le L
\quad\mbox{and}\quad
\max_\zs{1\le i\le d}
\sup_\zs{x,y\in\bbr}
\frac{|S_\zs{i}(y)-S_\zs{i}(x)|}{|y-x|}\le L\,.
$$
In this case (see, for example, \cite{Ga})  stochastic equation \eqref{A.9}
has the unique strong solution  $(y_\zs{t})_\zs{0\le t\le T}$ for any 
random variable $\theta$ with values in $\bbr^d$.

Moreover (see, for example \cite{LiSh}), for any $\theta\in\bbr^d$ the distribution
$\P_\zs{\theta}$ is absalutly continuous with respect to the Wiener measure $\nu_\zs{w}$
in $\cC[0,T]$ and the corresponding Radon-Nikodym derivative
for any function $x=(x_\zs{t})_\zs{0\le t\le T}$ from $\cC[0,T]$
 is defined as
\begin{equation}\label{A.11}
\frac{\d\P_\zs{\theta}}{\d \nu_\zs{w}}=
f(x,\theta)=\exp\left\{\int_0^T\,S(x_t,\theta)\d x_t -\frac{1}{2}\int_0^T\,S^2(x_t,\theta)\d t\right\}\,.
\end{equation}

Let $\Phi$ be a prior density in $\bbr^d$ having
the following form:
$$
\Phi(\theta)=\Phi(\theta_1,\ldots,\theta_d)=\prod_{j=1}^d\varphi_\zs{j}(\theta_j)\,,
$$
where $\varphi_\zs{j}$ is some continuously differentiable density in $\bbr$. 
Moreover, let $\lambda(\theta)$ be a continously differentiable $\bbr^d\to \bbr$ function such that
for each $1\le j\le d$
\begin{equation}\label{A.12}
\lim_\zs{|\theta_\zs{j}|\to\infty}\,
\lambda(\theta)\,\varphi_\zs{j}(\theta_\zs{j})=0
\quad\mbox{and}\quad
\int_\zs{\bbr^d}\,|\lambda^{\prime}_\zs{j}(\theta)|\,\Phi(\theta)\,\d \theta
<\infty\,,
\end{equation}
where
$$
\lambda^{\prime}_\zs{j}(\theta)=\frac{\partial\lambda(\theta)}{\partial\theta_j}\,.
$$
For any $\cB(\cX)\times\cB(\bbr^d)-$
measurable integrable function $\xi=\xi(x,\theta)$ we denote
$$
\tilde{\E}\xi=\int_{\bbr^d}\,
\int_\zs{\cX}\xi(x,\theta)\,\d \P_\zs{\theta}\,\Phi(\theta) \d \theta
=
\int_{\bbr^d}\,\int_\zs{\cX}\,
\xi(x,\theta)\,f(x,\theta)\,\Phi(\theta)\d \nu_\zs{w}(x)\, \d \theta\,,
$$
where $\cX=\cC[0,T]$.

\begin{lemma}\label{Le.A.2}
For any square integrable function $\hat{\lambda}_\zs{T}$
 measurable with respect to $(Y_t)_\zs{0\le t\le T}$
 and for any $1\le j\le d$
  the following inequality holds
$$
\tilde{\E}(\hat{\lambda}_\zs{T}-\lambda(\theta))^2\ge
\frac{\Lambda^2_\zs{j}}
{\tilde{\E}\int_0^T\,S^2_\zs{j}(Y_t)\,\d t +I_j}\,,
$$
where
$$
\Lambda_\zs{j}=\int_\zs{\bbr^d}\,\lambda^{\prime}_\zs{j}(\theta)\,\Phi(\theta)\,\d \theta
\quad\mbox{and}\quad
I_j=\int_\zs{\bbr}\,\frac{\dot{\varphi}^2_j(z)}{\varphi_\zs{j}(z)}\,\d z\,.
$$
\end{lemma}
{\bf  Proof.}
First of all note that for the function \eqref{A.10} and for the Wiener process
$w=(w_\zs{t})_\zs{0\le t\le T}$ density \eqref{A.11} is bounded 
with respect to $\theta_\zs{j}\in\bbr$ for any $1\le j\le d$, i.e.
$$
\limsup_\zs{|\theta_\zs{j}|\to\infty}\,f(w,\theta)\,<\,\infty
\quad\quad\mbox{a.s.}
$$
Therefore taking into account condition \eqref{A.12}
by integration by parts
one gets 
\begin{align*}
\tilde{\E}\left((\hat{\lambda}_\zs{T}-\lambda(\theta))\Psi_j\right)
&=\int_{\cX\times\bbr^d}\,(\hat{\lambda}_\zs{T}(x)-\lambda(\theta))\frac{\partial}{\partial\theta_j}
\left(f(x,\theta)\Phi(\theta)\right)\d \theta\nu_w(\d x)\\
&=\int_{\cX\times\bbr^{d-1}}\left(\int_{\bbr}\,
\lambda^{\prime}_\zs{j}(\theta)\,
f(x,\theta)\Phi(\theta)\d \theta_j\right)\prod_{i\neq j}\d \theta_i\nu_w(\d x)\\
&=\Lambda_\zs{j}\,.
\end{align*}
Now by the Bouniakovskii-Cauchy-Schwarz inequality
we obtain tha following lower bound for the quiadratic risk
$$
\tilde{\E}(\hat{\lambda}_\zs{T}-\lambda(\theta))^2\ge
\frac{\Lambda^2_\zs{j}}{\tilde{\E}\Psi_j^2}\,,
$$
where
\begin{align*}
\Psi_j=\Psi_\zs{j}(x,\theta)&=
 \frac{\partial}{\partial\theta_j}\,\ln(f(x,\theta)\Phi(\theta))\\
&=
 \frac{\partial}{\partial\theta_j}\,\ln f(x,\theta)
+
 \frac{\partial}{\partial\theta_j}\,\ln \Phi(\theta)\,.
\end{align*}
Note that 
from \eqref{A.11} it is easy to deduce 
that 
$$
\frac{\partial}{\partial\theta_j}\,\ln f(y,\theta)
=
\int^{T}_\zs{0}\,S_\zs{j}(y_\zs{t})\d w_\zs{t}\,.
$$
Therefore, due to the boundness of the functions $S_\zs{j}$
we find that for each $\theta\in\bbr^d$
$$
\E_\zs{\theta}\,
\frac{\partial}{\partial\theta_j}\,\ln f(y,\theta)\,=0
\quad\mbox{and}\quad
\E_\zs{\theta}\,
\left(
\frac{\partial}{\partial\theta_j}\,\ln f(y,\theta)\,
\right)^2
=\E_\zs{\theta}\,
\int^{T}_\zs{0}\,S^2_\zs{j}(y_\zs{t})\d t\,.
$$
Taking this into account we can calculate now 
$\tilde{\E}\Psi_j^2$, i.e.
$$
\tilde{\E}\Psi_j^2=
\tilde{\E}\,
\int^{T}_\zs{0}\,S^2_\zs{j}(y_\zs{t})\d t\,
+\,I_\zs{j}\,.
$$
Hence
Lemma~\ref{Le.A.2}.
\endproof

\subsection{Proof of (\ref{Lo.14})}\label{Su.A.5}

We set
$$
\psi_\zs{m,j}(y)=\frac{1}{h}D^2_\zs{m,j}(y)\,.
$$
Then by making use of definitions in \eqref{A.1}
we can estimate the term $\varsigma_\zs{m,j}(T)$ as
\begin{align*}
|\varsigma_\zs{m,j}(T)|&\le 
\frac{\E\left(
\E_\zs{S_\zs{\theta,T}}\,|\Delta_\zs{T}(\psi_\zs{m,j})|
\right)}{\ov{e}_\zs{j}(I^2_\zs{\eta})q_\zs{0}(\tilde{x}_\zs{m})\sqrt{T}}
+\E\,
\left|
\frac{M_\zs{S_\zs{\theta,T}}(\psi_\zs{m,j})}{\ov{e}_\zs{j}(I^2_\zs{\eta})q_\zs{0}(\tilde{x}_\zs{m})}
-1
\right|\,.
\end{align*}
Moreover, taking into account that
$$
\lim_\zs{\eta\to 0}\,\sup_\zs{j\ge 1}
|\ov{e}_\zs{j}(I^2_\zs{\eta})-1|=0
$$
we chose $\eta>0$ for which
$$
\inf_\zs{j\ge 1}\,\ov{e}_\zs{j}(I^2_\zs{\eta})\ge 1/2\,.
$$
Therefore we can write that
\begin{align}\nonumber
|\varsigma_\zs{m,j}(T)|&\le 
\frac{2}{q_\zs{*}\sqrt{T}}
\E\left(
\E_\zs{S_\zs{\theta,T}}\,|\Delta_\zs{T}(\psi_\zs{m,j})|
\right)
+
\frac{2}{q_\zs{*}}
\E\,
\left|
M_\zs{S_\zs{\theta,T}}(\psi_\zs{m,j})
-
M_\zs{S_\zs{0}}(\psi_\zs{m,j})
\right|\\ \label{A.13}
&+
\frac{2}{q_\zs{*}}
\left|
M_\zs{S_\zs{0}}(\psi_\zs{m,j})
-
\ov{e}_\zs{j}(I^2_\zs{\eta})q_\zs{0}(\tilde{x}_\zs{m})
\right|\,.
\end{align}
We remind that on the set 
\eqref{Lo.7} for sufficiently large $T$ the function 
$S_\zs{\theta,T}\in\Sigma_\zs{L}$
 therefore we estimatethe first term in the right side of the last inequality as
 $$
 \E\left(
\E_\zs{S_\zs{\theta,T}}\,|\Delta_\zs{T}(\psi_\zs{m,j})|
\right)\le \frac{2}{h}\,\P(\Xi^c_\zs{T})+
\sup_\zs{S\Sigma_\zs{L}}\,\E_\zs{S}
|\Delta_\zs{T}(\psi_\zs{m,j})|\,.
 $$
Moreover, taking into account that
$$
\int^{\infty}_\zs{-\infty}|\psi_\zs{m,j}(y)|\d y=
\int^1_\zs{-1}\,e^2_\zs{j}(v)\,I^2_\zs{\eta}(v)\d v\le 2
$$
we obtain that
$$
\limsup_\zs{T\to\infty}\max_\zs{1\le j\le M}\max_\zs{1\le j\le N}
\E\left(
\E_\zs{S_\zs{\theta,T}}\,|\Delta_\zs{T}(\psi_\zs{m,j})|
\right)\,<\,\infty\,.
$$
To estimate the next term in \eqref{A.13} we make use of the fact that on the set $\Xi_\zs{T}$
the function $S_\zs{\theta,T}$ satisfies inequality \eqref{Lo.9}
and one can check directly that on this set
$$
|q_\zs{S_\zs{\theta,T}}-q_\zs{0}|_\zs{*}
\le C^*\,(e^{2(b-a)\epsilon_\zs{T}}-1)\,.
$$
Therefore, with the help of this inequality we obtain
 that\begin{align*}
\E\,
\left|
M_\zs{S_\zs{\theta,T}}(\psi_\zs{m,j})
-
M_\zs{S_\zs{0}}(\psi_\zs{m,j})
\right|
&\le 
\int^1_\zs{-1}\,e^2_\zs{j}(v)I^2_\zs{\eta}(v)\,
\E\,|q_\zs{S_\zs{\theta,T}}(\tilde{x}_\zs{m}+vh)-q_\zs{0}(\tilde{x}_\zs{m}+vh)|\d v\\
&
\le \,2\P(\Xi^c_\zs{T})q^*\,+
2\,C^*\,(e^{2(b-a)\epsilon_\zs{T}}-1)\,.
\end{align*}
Finally it is easy to see that
\begin{align*}
\left|
M_\zs{S_\zs{0}}(\psi_\zs{m,j})
-
\ov{e}_\zs{j}(I^2_\zs{\eta})q_\zs{0}(\tilde{x}_\zs{m})
\right|\,&\le 
\int^1_\zs{-1}\,e^2_\zs{j}(v)I^2_\zs{\eta}(v)\,
|q_\zs{0}(\tilde{x}_\zs{m}+vh)-q_\zs{0}(\tilde{x}_\zs{m})|\d v\\
&
\le 2\sup_\zs{a\le u,v\le b\,,|u-v|\le h}\,
|q_\zs{0}(u)-q_\zs{0}(v)|\,,
\end{align*}
i.e. this term goes to zero as $h\to 0$ uniformly over $1\le m\le M$ and $1\le j\le N$.
Hence limit equality \eqref{Lo.14}

\endproof

\medskip

%\newpage

\begin{flushright}
\begin{tabular}{lcl}
   L.Galtchouk                       &$\quad$& S. Pergamenshchikov              \\
 IRMA, Department of Mathematics           &$\quad$& Laboratoire de Math\'ematiques,\\      
 Strasbourg University               &$\quad$& Avenue de l'Universit\'e, BP. 12,            \\
 7, st. Rene Descartes               &$\quad$&  Universit\'e de Rouen,                  \\
 67084, Strasbourg                   &$\quad$&  F76801, Saint Etienne du Rouvray\\
 Cedex, France                       &$\quad$&           Cedex, France\\
 e-mail: galtchou@math.u-strasbg.fr  &$\quad$& Serge.Pergamenchtchikov@univ-rouen.fr         \\
\end{tabular}
\end{flushright} 

\end{document}